\input amstex
\documentstyle{amsppt}
\magnification=\magstep1
\hsize=5.2in
\vsize=6.7in

\catcode`\@=11
\loadmathfont{rsfs}
\def\mycal{\mathfont@\rsfs}
\csname rsfs \endcsname

\topmatter
\title  A UNIQUE DECOMPOSITION RESULT \\
FOR HT FACTORS WITH TORSION FREE CORE
\endtitle
\author SORIN POPA \endauthor

\rightheadtext{Unique decomposition result}

\affil University of California, Los Angeles\endaffil

\address Math.Dept., UCLA, LA, CA 90095-155505\endaddress
\email popa@math.ucla.edu\endemail

\thanks Supported in part by NSF Grant 0100883.\endthanks

\abstract We prove that II$_1$ factors $M$ have a unique (up to
unitary conjugacy) cross-product type decomposition around ``core
subfactors'' $N \subset M$ satisfying the property HT of ([P1])
and a certain ``torsion freeness'' condition. In particular, this
shows that isomorphism of factors of the form $L_\alpha(\Bbb Z^2)
\rtimes \Gamma$, for torsion free, non-amenable subgroups
$\Gamma\subset SL(2, \Bbb Z)$ and $\alpha = e^{2\pi i t}$, $t
\not\in \Bbb Q$, implies isomorphism of the corresponding groups
$\Gamma$.

\endabstract

\endtopmatter

\document

Let $M$ be a type II$_1$ factor and $B \subset M$ a von Neumann
subalgebra. Recall from (Section 2 in [P1]) that $M$ has the {\it
property} H {\it relative to} $B$ if the identity map on $M$ can
be approximated pointwise by regular (i.e., subunital, subtracial)
completely positive maps on $M$ that ``vanish at infinity''
relative to $B$. Also, $B \subset M$ is a {\it rigid inclusion}
(or $(M, B)$ has the {\it relative property} (T)) if regular
completely positive maps on $M$ that are close (pointwise) to the
identity on $M$ follow uniformly close to the identity on $B$ (cf.
4.2 in [P1]). If $B \subset M$ satisfies both these conditions,
then it is called a {\it HT inclusion} and $B$ is called a {\it HT
subalgebra} of $M$. (N.B. In fact, the notation used in [P1] to
designate this property is ``HT$_s$'', while ``HT'' is used for a
slightly weaker condition. We opted for the notation ``HT'' in
this paper for the benefit of simplicity.) We refer to ([P1]) for
the detailed definitions, as well as for notations and terminology
used hereafter.

In ([P1]), one primarily studies HT Cartan subalgebras of $M$,
i.e., HT inclusions $B\subset M$ with $B$ maximal abelian in $M$
and with the normalizer $\Cal N_M(B) = \{u\in \Cal U(M)\mid uBu^*
= B\}$ generating $M$. Thus, the main technical result in ([P1])
shows the uniqueness, up to unitary conjugacy, of the HT Cartan
subalgebras in separable II$_1$ factors. Such unique decomposition
results are crucial for the calculation of invariants of II$_1$
factors (e.g., the fundamental and automorphism groups) and
ultimately for their classification. Thus, one application of the
uniqueness result for HT Cartan subalgebras in ([P1]) shows that
two cocycle group von Neumann factors of the form
$L_{\mu_{\alpha_j}}(\Bbb Z^2 \rtimes SL(2, \Bbb Z))$, with
$\mu_{\alpha_j}$ the $SL(2, \Bbb Z)$-invariant cocycle on $\Bbb
Z^2$ given by ``rational rotation'' $\alpha_j = e^{2 \pi i t_j},
t_j\in \Bbb Q$, $j=1,2$, are isomorphic if and only if the
rational numbers $t_1, t_2$ have the same denominator.

In this paper we prove a new
unique decomposition result, this time for HT inclusions $B \subset M$ with
factorial  ``core'' subalgebras, $B=N$, which satisfy the following
``strong irreducibility'' condition:

\vskip .1in \noindent {\bf Definition}. An inclusion $N \subset M$
is {\it torsion free} if $N$ is a factor and $Q'\cap M= \Bbb C$
for any subfactor $Q \subset N$ with $[N:Q] < \infty$ and $Q'\cap
N = \Bbb C$. By a Proposition below, an inclusion of the form $N
\subset M=N\rtimes_\sigma \Gamma$ is torsion free iff $N$ is a
factor, the action $\sigma$ is properly outer and the group
$\Gamma$ is torsion free, a fact that justifies the terminology.

If $N \subset M$ is torsion free and HT, then we also say that
$N$ is a {\it torsion free} HT {\it core} for $M$.
We denote by $\Cal H\Cal T_{_{tf}}$ the class of factors that admit
torsion free HT cores.
The following are
examples of such (inclusions of) factors:

\vskip .1in \noindent {\bf Example 1.} Let $M_\alpha(\Gamma)$
denote the factors defined as in (3.3.2 and 6.9.1 of [P1]), i.e.,
$M_\alpha(\Gamma) =R_\alpha \rtimes_{\sigma_\alpha} \Gamma$, where
$\alpha = e^{2\pi i t}$ for some $t \in [0, 1/2]\setminus \Bbb Q$,
$R_\alpha$ is the hyperfinite II$_1$ factor represented as the
``irrational rotation algebra'' $L_\alpha(\Bbb Z^2)$, $\Gamma
\subset SL(2, \Bbb Z)$ is a non-amenable (equivalently
non-solvable) group and $\sigma_\alpha$ is the action of $\Gamma$
on $L_\alpha(\Bbb Z^2)$ induced by the action of $SL(2, \Bbb Z)$
on $\Bbb Z^2$. Note that by (page 62 of [Bu]) a group $\Gamma
\subset SL(2, \Bbb Z)$ is non-amenable if and only if the pair
$(\Bbb Z^2 \rtimes \Gamma, \Bbb Z^2)$ has the relative property
(T) ([M]). Thus, by (5.1 in [P1]) the inclusions $R_\alpha \subset
M_\alpha(\Gamma)$ are rigid and since all $\Gamma \subset SL(2,
\Bbb Z)$ have Haagerup's compact approximation property, by (3.1
in [P1]), $M_\alpha(\Gamma)$ has property H relative to
$R_\alpha$. Moreover, by (3.3.2 $(ii)$ in \cite{Betti})
$\sigma_\alpha$ are properly outer actions (because the action of
$\Gamma$ on $\Bbb Z^2$ is outer) and by (3.3.2 $(i)$ in
\cite{Betti}) they are ergodic as well. Indeed, this is because
the stabilizer of any non-zero element in $\Bbb Z^2$ is a cyclic
subgroup of $SL(2, \Bbb Z)$, so if $\Gamma$ leaves a finite subset
$\neq \{(0,0)\}$ of $\Bbb Z^2$ invariant, then it is almost
cyclic.

Altogether, this shows that $R_\alpha \subset M_\alpha(\Gamma)$
are irreducible HT inclusions of II$_1$ factors. The algebras
$M_\alpha(\Gamma)$, will be called {\it irrational rotation} HT
{\it factors}. Note that the inclusion $R_\alpha \subset
M_{\alpha}(\Gamma)$ is torsion free whenever the non-amenable
group $\Gamma\subset SL(2, \Bbb Z)$ is torsion free. In particular
one can take $\Gamma=\Bbb F_n$ for any $ 2 \leq n\leq \infty$.

More generally, if $\Gamma$ is a torsion free Haagerup group
acting outerly on a group $H$ such that the pair $(H \rtimes
\Gamma, H)$ has the relative property (T), and if $\mu$ is a
$\Gamma$-invariant scalar 2-cocycle on $H$ such that $L_\mu(H)$ is
a factor, then $L_\mu(H) \rtimes \Gamma \in \Cal H\Cal T_{_{tf}}$,
with $L_\mu(H)$ a torsion free HT core for $L_\mu(H) \rtimes
\Gamma$. The irrational rotation factors $M_\alpha(\Gamma)$
correspond to the case the cocycle $\mu=\mu_\alpha$ is given by
$\alpha=e^{2\pi i t}$, $t\not\in \Bbb Q$.

\vskip .1in
\noindent
{\bf Example 2}. Let $N$ be a type II$_1$
factor with the property (T) of Connes-Jones ([C1,2], [CJ]). If
$\beta$ is an aperiodic automorphism of $N$ and one still denotes
by $\beta$ the properly outer action of $\Bbb Z$ on $N$
implemented by $\{\beta^n\}_{n \in \Bbb Z}$, then $N \subset N
\rtimes_\beta \Bbb Z$ is a torsion free HT inclusion. More
generally, any properly outer cocycle action $\sigma : \Gamma
\rightarrow {\text{\rm Aut}}(N)$ of a torsion free group $\Gamma$
satisfying Haagerup's compact approximation property gives rise to
a torsion free HT inclusion $N \subset N\rtimes_\sigma \Gamma$.
Indeed, by (4.7 and 5.9 in [P1]) $N \subset M$ is rigid, while $M$
has the property H relative to $N$ by (3.1 in [P1]). Note that
these examples include some of the factors considered in the last
part of Example 1: Thus, if $H$ is taken to be an ICC group with
the property (T) of Kazhdan and $\Gamma$ a torsion free Haagerup
group of outer automorphisms of $H$ (e.g., $\Gamma=\Bbb Z$), then
$N=L(H) \subset L(H\rtimes \Gamma) =M$ is a torsion free HT
inclusion. More generally, one can take $H$ arbitrary with the
property (T), but with a scalar $\Gamma$-invariant cocycle $\mu$
on it such that $N=L_\mu(H)$ is a factor.

\vskip .1in
\noindent
{\bf Example 3}. If $M$ is a II$_1$ factor
with torsion free HT core $N \subset M$ and $t > 0$ then the
amplification $M^t$ of $M$ has $N^t$ as HT core (by 2.4 and 4.7 in
[P1]), which trivially follows torsion free. Also, if $N_i \subset
M_i = N_i \rtimes \Gamma_i, i=1,2,$ are cross-product torsion free
HT inclusions then $N_1\overline{\otimes} N_2 \subset M_1
\overline{\otimes} M_2$ is torsion free HT. Thus, the class of
factors with torsion free HT core is well behaved to
amplifications and tensor products, producing more examples from
the ones described in 1 and 2 above.

The main result in this paper shows that II$_1$ factors in the class
$\Cal H \Cal T_{_{tf}}$ have unique decomposition
(up to unitary conjugacy) around their
torsion free HT cores:

\proclaim{Theorem} If $N_1, N_2 \subset M$ are torsion free
${\text{\rm HT}}$ inclusions of $\text{\rm II}_1$ factors
then there exists a unitary
element $u\in M$ such that $uN_1u^*=N_2$. Thus,
factors in the class $\Cal H\Cal T_{_{tf}}$ have unique (up to unitary conjugacy)
torsion free ${\text{\rm HT}}$ core.
\endproclaim

The case of interest is when $M\in \Cal H\Cal T_{_{tf}}$
are cross-product factors
over their torsion free HT cores,
when the above statement becomes:

\proclaim{Corollary 1} Let $N_i$ be $\text{\rm II}_1$ factors,
$\Gamma_i$  be torsion free groups with Haagerup property
and $\sigma_i$ be properly outer cocycle actions of $\Gamma_i$
on $N_i$ such that $N_i \subset N_i \rtimes_{\sigma_i}
\Gamma_i$ are rigid inclusions, $i=1,2$. If
$N_1 \rtimes_{\sigma_1}
\Gamma_1 \simeq N_2 \rtimes_{\sigma_2} \Gamma_2$, via some
isomorphism $\theta$,
then there exists a unitary element
$u \in N_2 \rtimes_{\sigma_2} \Gamma_2$
such that ${\text{\rm Ad}}(u)\circ \theta$ takes
$N_1$ onto $N_2$, implements an isomorphism between $\Gamma_1$,
$\Gamma_2$ and cocycle conjugates
$\sigma_1, \sigma_2$.
\endproclaim

In particular, Corollary 1 shows that non-isomorphic torsion-free
groups $\Gamma$ give rise to non-isomorphic factors $N \rtimes
\Gamma$ in the class $\Cal H\Cal T_{_{tf}}$. For the irrational
rotation HT factors in Example 1, this gives: If $\Gamma_1,
\Gamma_2$ are torsion free non-amenable subgroups of $SL(2, \Bbb
Z)$ and $\Gamma_1 \not\simeq \Gamma_2$ then
$M_{\alpha_1}(\Gamma_1) \not\simeq M_{\alpha_2}(\Gamma_2)$,
$\forall \alpha_1, \alpha_2$. The Corollary also shows that if
$N_i$ are property (T) II$_1$ factors with aperiodic automorphisms
$\beta_i, i=1,2$, then $N_1\rtimes_{\beta_1} \Bbb Z \simeq
N_2\rtimes_{\beta_2} \Bbb Z$ iff $N_1 \simeq N_2$ and $\beta_1$
cocycle conjugate to $\beta_2$ (cf. Example 2).

By the Theorem above, classical invariants for the
factors $M\in \Cal H\Cal T_{_{tf}}$, such as the automorphism group
Out$(M)$ or the fundamental
group $\mycal F(M)$, coincide with their ``relative'' versions,
Out$(N \subset M)={\text{\rm Aut}}(N\subset M)/\{\text{\rm Ad}(u) \mid
u \in \Cal N_M(N)\}$ and respectively $\mycal F(N\subset M)
= \{t>0 \mid (N\subset M)^t \simeq (N\subset M)\}$,
whenever $N \subset M$ is a torsion free HT core for $M$.

Moreover,
if we denote by $\Cal G_N$ the group
of automorphisms of $M$
generated by $\text{\rm Int}(M)$ and by the automorphisms
that leave $N$ pointwise fixed
($\simeq \hat{\Gamma}$ when $M= N \rtimes \Gamma$, by [JP]),
then by (4.4 of [P1]) $\Cal G_N$ is open and closed
in Aut$(M)$ so
the quotient group Aut$(M)/\Cal G_N$ is countable. By the Theorem,
this quotient group
is an isomorphism invariant
for the factors in the class $\Cal H\Cal T_{_{tf}}$.
We denote this invariant by
Out$_{_{HT}}(M)$, for $M \in \Cal H\Cal T_{_{tf}}$ (same notation as for
its analogue invariant for the class $\Cal H\Cal T$ in [P1]).

Since $N\subset M$ torsion free HT implies
$N\overline{\otimes} N \subset M\overline{\otimes} M$
torsion free HT, and since a choice of $\theta^t : M \simeq M^t$
for each $t \in \mycal F(M)$ gives a one to one embedding
$\{ \theta^t \otimes \theta^{1/t}\mid t \in \mycal F(M)\}\subset$
Out$_{_{HT}}(M\overline{\otimes} M)$
(cf. [C1]), it follows
that $\mycal F(M)$ is countable as well.

The Theorem shows that if
$N \subset M=N\rtimes_\sigma \Gamma$ is a
torsion free HT inclusion then any automorphism $\theta$
of $M$ (respectively of $M^\infty \overset \text{\rm def} \to =
M\overline{\otimes} \Cal B(\ell^2\Bbb N)$),
can be perturbed by an automorphism in $\Cal G_N$ (resp.
in the group of automorphisms of $M^\infty$ that are amplifications of
automorphisms in $\Cal G_N$) to
an automorphism $\theta'$ that takes the core $N$
(resp. $N^\infty$) onto itself.
Thus, $\theta'_{|N}$ (resp. $\theta'_{|N^\infty}$)
must lie in the normalizer
$\Cal N(\sigma)$ of $\sigma(\Gamma)$ in Out$(N)$ (resp. in the
normalizer $\Cal N^\infty(\sigma)$ of $\sigma(\Gamma)$
in Out$(N^\infty)$). Since conversely any automorphism in
$\Cal N$ implements an automorphism
of $N \rtimes_\sigma \Gamma$, it follows that
Out$_{HT}(M)$ is isomorphic
to $\Cal N(\sigma)/\sigma(\Gamma)$ and
$\Cal F(M)$ is isomorphic to the image
via Mod of $\Cal N^\infty(\sigma)$,
i.e., $\Cal F(M)=$Mod$(\Cal N^\infty(\sigma))$.

Moreover, Corollary 1 implies that
$\Cal C(\sigma)\overset \text{\rm def} \to =
\sigma(\Gamma)'
\cap {\text{\rm Out}}(N)$ and the group $\Cal O(\sigma)$
of outer automorphisms of the group $\Gamma$ implemented by
elements in $\Cal N(\sigma)$ are isomorphism invariants
for $N\rtimes_\sigma \Gamma \in \Cal H\Cal T_{_{tf}}$.
Altogether, we have:

\proclaim{Corollary 2} If $M$ is a $\text{\rm II}_1$ factor with
torsion free $\text{\rm HT}$ core $N \subset M$ then
$\mycal F(M)=\mycal F(N\subset M)$, $\text{\rm Out}(M)=
\text{\rm Out}(N \subset M)$. Moreover, $\mycal F(M)$
and ${\text{\rm Out}}_{_{HT}}(M)$ are countable. If
in addition $M = N \rtimes_\sigma \Gamma$, then
$\text{\rm Out}_{_{HT}}=\Cal N(\sigma)/\sigma(\Gamma)$,
$\mycal F(M)={\text{\rm Mod}}(\Cal N^\infty(\sigma))$. Also,
$\Cal C(\sigma), \Cal O(\sigma)$ are isomorphism
invariants for $M$ and are countable.
\endproclaim

To prove the Theorem, we need two lemmas. The first one
shows that in order for the irreducible,
torsion free subfactors $N_1, N_2 \subset M$
to be conjugate in $M$, it is sufficient
to have a finite dimensional ``intertwining'' bimodule
between them.

\proclaim{Lemma 1} Let $N_1, N_2 \subset M$ be irreducible subfactors
(i.e., $N_i'\cap M = \Bbb C, i=1,2$).
Assume $N_2 \subset M$
is torsion free. If there exists a non-zero,
finite dimensional  $N_1-N_2$ sub-bimodule of $L^2(M)$,
then there exists $u \in \Cal U(M)$ such that $uN_1u^* \subset N_2$
with $[N_2: uN_1u^*] < \infty$. If in addition $N_1\subset M$
is torsion free as well, then any $u$ as above must satisfy $uN_1u^*=N_2$.
\endproclaim
\noindent
{\it Proof}. This can be easily derived from (2.1 in [P2]), but we'll
give here a self contained
argument. Let $\Cal H \subset L^2(M)$ be so that
dim$_{N_1}\Cal H$, dim$\Cal H_{N_2} < \infty$. By taking a submodule
of $\Cal H$ if necessary, we may assume $\Cal H$ is irreducible.
Thus, $N_1 \subset JN_2J' \cap \Cal B(\Cal H)$
is an irreducible inclusion of finite index. Equivalently,
if we denote by $p$ the projection
of $L^2(M)$ onto $\Cal H$ then $p \in N_1' \cap
\langle M, N_2 \rangle$, $Tr(p) < \infty$ and the
inclusion $N_1p \subset
p\langle M, N_2 \rangle p$ is irreducible with finite index.

Since $N_1$ is of type II$_1$, there exists $q_1\in \Cal P(N_1),
q_1\neq 0,$ such that
$Tr(pq_1) \leq 1$. Since $Tr(e_{N_2})=1$, it follows
that $pq_1$ is majorized by $e_{N_2}$ in the
type II factor $\langle M,
N_2 \rangle$. But $e_{N_2}\langle M,
N_2 \rangle e_{N_2}=N_2e_{N_2}$, so there exists $q_2 \in N_2$ and
a partial isometry $V \in \langle M, N_2 \rangle$ such that $V^*V=pq_1$,
$VV^* = e_{N_2}q_2$. By spatiality,  $V(q_1N_1q_1p)V^*$ is an
irreducible subfactor of finite index in $q_2N_2q_2e_{N_2}$.

Let $Q \subset N_2$ be such that $q_2 \in Q$ and $q_2Qq_2e_{N_2}
= V(q_1N_1q_1p)V^*$ and denote by $\theta: q_1N_1q_1
\simeq q_2Qq_2$ the isomorphism satisfying $VxV^* = \theta(x)e_{N_2}$.
Equivalently, $Vx = \theta(x)V, \forall x\in q_1N_1q_1$.
By applying the canonical
operator valued weight $\Phi$ of $\langle M, N_2 \rangle$ onto $M$,
it follows that $\xi = \Phi(V)$, which apriorically
lies in $L^2(M)$, satisfies $e_{N_2}\xi  =V$ and
$\xi x = \theta(x) \xi, \forall x\in q_1N_1q_1$. Thus $\xi^*\xi x =
\xi^* \theta(x) \xi = x \xi^* \xi, \forall x$, and since
$(q_1N_1q_1)'\cap q_1Mq_1 = \Bbb Cq_1$, it follows that $\xi$
is a scalar multiple of a partial isometry $v \in M$ with
$v^*v=q_1$.

Similarly, the intertwiner relation also gives
$vv^* \in (q_2Qq_2)'\cap q_2Mq_2$. But
by the torsion freeness of $N_2 \subset M$ we have
$(q_2Qq_2)'\cap q_2Mq_2=\Bbb Cq_2$. Thus, $v(q_1N_1q_1)v^*
=q_2Qq_2 \subset q_2N_2q_2$ and since both $N_1, N_2$ are factors,
there exists a unitary element $u\in M$ with $uq_1=v$ and
$uN_1u^* \subset N_2$.

Finally, if $N_1 \subset M$ is torsion free as well, then
let $Q_1 \subset N_1$ be a ``downward basic
construction'' for $N_1 \subset u^*N_2u$ (cf. [J]). It follows that
$Q_1'\cap N_1 = \Bbb C$ but $Q_1'\cap M \supset Q_1'\cap u^*N_2u \neq \Bbb C$
unless $N_1=u^*N_2u$.
\hfill $\square$

\proclaim{Lemma 2} Let $M$ be a separable ${\text{\rm II}}_1$ factor.
Assume $B_1, B_2 \subset M$ are von Neumann subalgebras
such that $M$ has the property ${\text{\rm H}}$
relative to $B_1$ and $B_2 \subset M$ is rigid. Then
$B_2$ is discrete over $B_1$, i.e.,
$L^2(M)$ is generated by $B_2-B_1$ bimodules which are
finite dimensional over $B_1$.
\endproclaim
\noindent
{\it Proof}. By the property H of $M$ relative to $B_1$ there
exist regular,
completely positive, $B_1$-bimodular
maps  $\phi_n$ on $M$ such that $\phi_n \rightarrow id_M$ and
$T_{\phi_n} \in \Cal J_0(\langle M, B_1 \rangle)$. By the rigidity
of $B_2\subset M$ it follows that $\varepsilon_n=$
sup$\{\|\phi_n(u) - u\|_2 \mid u\in \Cal U(B_2)\} \rightarrow 0$.
Fix $x\in M$ and note that by (Corollary 1.1.2 of [P1])
we have
$$
\|u^*T_{\phi_n}u(\hat{x})-\hat{x}\|_2=\|\phi_n(ux) -ux\|_2
$$
$$
\leq \|\phi_n(ux) -u\phi_n(x)\|_2 + \|\phi_n(x)-x\|_2 \leq
2\varepsilon_n^{1/2} + \|\phi_n(x)-x\|_2.
$$
Thus, by taking weak limits of appropriate convex combinations
of elements
of the form $u^*T_{\phi_n}u$ with $u\in \Cal U(B_2),$ and
using (Proposition 1.3.2 of [P1]), it follows that
$T_n=\Cal E_{B_2'\cap \langle M, B_1\rangle} (T_{\phi_n})
\in K_{T_{\phi_n}} \cap
(B_2'\cap J_0(\langle M, B_1 \rangle))$ satisfy
$\underset n \rightarrow \infty \to \lim
\|T_n(\hat{x})-\hat{x}\|_2 = 0$. But $x\in M$ was arbitrary.
This shows that the right supports of $T_n$ span all the identity
of $\langle M, B_1 \rangle$. Since $T_n$ are compact,
this shows that $B_2'\cap \langle M, B_1\rangle$
is generated by finite projections of $\langle M, B_1\rangle$.
Equivalently, $B_2$ is discrete relative to $B_1$.
\hfill $\square$

\vskip .1in
\noindent
{\it Proof of the Theorem}. By Lemma 2, $N_1$ is
discrete over $N_2$ and $N_2$ is discrete over $N_1$.
Thus, $L^2(M)$ is generated by finite
dimensional $N_1-N_2$ bimodules. By the torsion freeness
of $N_1, N_2 \subset M$ and Lemma 1, this implies
$N_1, N_2$ are unitary conjugate.
\hfill $\square$

\vskip .1in
We'll now discuss in more details the
torsion freeness condition. In particular, we
prove that in the case of
cross product inclusions, this condition
amounts to the group involved being torsion free.

First recall some terminology
from ([P3,1]):
The {\it quasi-normalizer}
of a subfactor $N \subset M$
is the set $q\Cal N_M(N)=\{x\in M
\mid \text{\rm dim}_NL^2(NxN)_N < \infty\}$. Note that the linear
span of $q\Cal N_M(N)$ is a $*$-subalgebra of $M$
containing $N$. Also, note that if $Q \subset N$
is a subfactor with $[N:Q] < \infty$ then $q\Cal N_M(Q)=
q\Cal N_M(N)$ (see [P1]). $N$ is {\it quasi-regular} (or
{\it discrete}) in $M$ if $q\Cal N_M(N)''=M$.

A typical example of quasi-regular subalgebras is
when $N'\cap M = \Bbb C$ and $N$ is {\it regular} in $M$,
i.e., $\Cal N_M(N)''=M$, equivalently when
$M = N \rtimes_\sigma \Gamma$ for some properly outer cocycle action
of the group $\Gamma=\Cal N(N)/\Cal U(N)$ on $N$. Other examples
are the symmetric enveloping inclusions associated to
extremal subfactors of finite Jones index ([P3]).

Note that by
(3.4 in ([P1]), if $M$ has the property H relative to
a subfactor $N\subset M$  then $N$ is quasi-regular
in $M$. Thus, a HT inclusion of factors $N \subset M$
as in the Theorem is automatically quasi-regular.

\proclaim{Proposition} $1^\circ$. Let $N\subset M$ be an irreducible
inclusion of factors. If $q\Cal N_M(N)''=\Cal N_M(N)''$ and we
denote $\Gamma=\Cal N_M(N)/\Cal U(N)$, then
$N \subset M$ is torsion free iff $\Gamma$ is torsion free.

$2^\circ$. If $N$ is quasi-regular in $M$ and
$N\subset M$ is torsion free then there
exist no intermediate subfactors $N\subset P \subset M$ such that
$[P:N] < \infty$, $P\neq N$.
\endproclaim
\noindent
{\it Proof}.
Indeed, for if $u_0\in \Cal N_M(N)$ implements an
automorphism $\theta_0$ with outer period $2\leq k < \infty$ then
by Connes' theorem there exists $v \in \Cal U(N)$ such that $(vu_0)^k=1$.
Thus, if $\theta=$Ad$(vu_0)$ then $Q=N^\theta$ satisfies
$Q'\cap N=\Bbb C$ while $vu_0 \in (Q'\cap M) \setminus \Bbb C1$.

Conversely, let $Q$ be an
irreducible subfactor of $N$ with finite index. Assume
$a \in (Q'\cap M) \setminus \Bbb C1$. Then $a \in q\Cal N_M(Q)
= q\Cal N_M(N)$ (the equality
holds because $Q$ has finite index in $N$).
By hypothesis, it follows that $a = \Sigma_g x_gu_g $ for some
$x_g \in N$, where $u_g \in \Cal N_M(N)$ are some
unitaries implementing the cross-product
construction $N \subset N \rtimes \Gamma = \Cal N_M(N)''$,
where $\Gamma = \Cal N_M(N)/\Cal U(N)$. This implies that
there exists $g \neq e$ such that $x_g \neq 0$ and
$yx_g = x_g \sigma_g(y), \forall y\in Q$. By $Q'\cap N=\Bbb C$ it follows
that the partial isometry in the polar decomposition of $x_g$
is a unitary element $v$ with $\theta=$Ad$(vu_g)$ satisfying $Q \subset N^\theta$.
Since $[N:Q] < \infty$, this shows that $\theta$ is periodic, thus $g$ has torsion.

$2^\circ$. If there exists $N\subset P \subset M$ with
$[P:N] < \infty$ and $P\neq N$ then let
$Q \subset N$ be a downward basic construction for $N \subset P$
(cf. [J]). We then have $Q'\cap N = \Bbb C$ (because
$N'\cap P = N'\cap M = \Bbb C$) but $Q'\cap P \neq \Bbb C$, as
it contains the Jones projection.
\hfill $\square$

\vskip .1in
\noindent
{\bf Remarks. 1}$^\circ$. The converse implication in part 2$^\circ$ of
the Proposition is probably true as well. This is of course the case
when $M = N \rtimes_\sigma \Gamma$, by part 1$^\circ$ of that statement.

{\bf 2$^\circ$}. By Corollary 1, isomorphism of
cross product HT factors $M=N \rtimes_\sigma \Gamma$
(such as the irrational rotation HT factors $M_\alpha(\Gamma)$), with
torsion free $\Gamma$, amounts to
the isomorphism of the groups $\Gamma$ and the cocycle
conjugacy of the (cocycle) actions $\sigma$.
But in both Examples 1 and 2,
we could not find ways to completely distinguish between the
(cocycle) actions $\sigma$.
In particular, we could  not obtain
precise calculations of
Out$(M)$, $\mycal F(M)$ (or for that matter
$\Cal C(\sigma)$, $\Cal O(\sigma)$)
by the method of calculation
of Corollary 2. However,
a classification (= non-isomorphism) ``modulo countable
sets'' of the factors $M_\alpha(\Gamma)$
is obtained in ([NPS]). Note that these factors (and in fact
all non-McDuff factors
of the form $N \rtimes_\sigma \Gamma$ with $N \simeq R$
and $\Gamma$ in the class $\Cal C$ of [O])
follow prime by Ozawa's recent results in ([O]).

{\bf 3$^\circ$}. Related to Example 2,
one can prove the following statement, by using an argument similar to
the proof of (Proposition 9 in [GP]):
If $(N, \tau)$ is a finite von Neumann
algebra with a faithful normal trace state and $\sigma$ is a
properly outer cocycle action of an amenable group $\Gamma$
on $(N, \tau)$, then the inclusion $(N \subset N\rtimes_\sigma \Gamma)$
is rigid if and only if $N$ has the property $(T)$ in ([P1]),
i.e., iff $N \subset N$ is rigid. This notion of property (T)
coincides with the one considered in ([Jol]) in the case the algebras are
of the form $N=L_\mu(H)$, for $H$ a discrete group
and $\mu$ a cocycle on it, when in fact
both conditions are
equivalent to the property (T) for $H$.
They also coincide in the case $N$ has finite dimensional center,
when they are equivalent to
the original Connes-Jones definition in ([CJ]).

{\bf 4$^\circ$}. It is interesting to know whether the free group
factors $L(\Bbb F_n)$ can be realized as ``cores'' of HT (or
merely rigid) inclusions $L(\Bbb F_n) \subset M$. It may be that
this happens iff $n$ is finite. The affirmative answer to the
``if'' part of this problem alone (i.e., showing that $L(\Bbb
F_2)$ can be realized as a rigid core), could provide new insight
to the ``(non)isomorphism of the free group factors'' problem. In
this respect, we should mention that it is an open problem whether
$\Bbb F_n \subset {\text{\rm Aut}}(\Bbb F_n)$ has the relative
property (T) ([M]) for some $2 \leq n < \infty$. In particular, it
is not known whether $\Bbb F_2 \subset \Bbb F_2 \rtimes SL(2, \Bbb
Z)$ has the relative property (T) or not.

\head  References\endhead

\item{[Bu]} M. Burger, {\it Kazhdan constants for} $SL(3,\Bbb Z)$,
J. reine angew. Math., {\bf 413} (1991), 36-67.

\item{[C1]} A. Connes: {\it A type II$_1$
factor with countable fundamental group}, J. Operator
Theory {\bf 4} (1980), 151-153.

\item{[C2]} A. Connes: {\it Classification des facteurs},
Proc. Symp. Pure Math.
{\bf 38}
(Amer. Math. Soc., 1982), 43-109.

\item{[CJ]} A. Connes, V.F.R. Jones: {\it Property} (T)
{\it for von Neumann algebras}, Bull. London Math. Soc. {\bf 17} (1985),
57-62.

\item{[GP]} D. Gaboriau, S. Popa: {\it An Uncountable Family of
Non Orbit Equivalent Actions of $\Bbb F_n$}, preprint,
math.GR/0306011.

\item{[dHV]} P. de la Harpe, A. Valette: ``La propri\'et\'e T
de Kazhdan pour les
groupes localement compacts'', Ast\'erisque {\bf 175}, Soc. Math. de France (1989).

\item{[Jol]} P. Jolissaint: {\it Property} (T) {\it for
discrete groups in terms of their regular representation}.

\item{[J]} V.F.R. Jones : {\it Index for subfactors}, Invent. Math.
{\bf 72} (1983), 1-25.

\item{[JP]} V.F.R. Jones, S. Popa: {\it Some properties of
MASAs in factors}, in ``Invariant subspaces and other topics'', pp. 89-102, Operator Theory: Adv. Appl. {\bf 6}, Birkhäuser, 1982.

\item{[K]} D. Kazhdan: {\it Connection of the dual space of a group
with the structure of its closed subgroups}, Funct. Anal. and its Appl.
{\bf1} (1967), 63-65.

\item{[M]} G. Margulis: {\it Finitely-additive invariant measures
on Euclidian spaces}, Ergodic. Th. and Dynam. Sys. {\bf 2} (1982),
383-396.

\item{[MvN]} F. Murray, J. von Neumann: {\it Rings of operators IV},
Ann. Math. {\bf44} (1943), 716-808.

\item{[NPS]} R. Nicoara, S. Popa, R. Sasyk :
{\it Some remarks  on irrational rotation $HT$ factors}, preprint,
math.OA/0401...

\item{[O]} N. Ozawa: {\it A Kurosh type theorem for} II$_1$ {\it
factors}, preprint, math.OA/0401121.

\item{[P1]} S. Popa: {\it On a class of type} II$_1$
{\it factors with Betti numbers invariants}, preprint OA/0209130.

\item{[P2]} S. Popa: {\it
Strong Rigidity of} II$_1$
{\it Factors Coming from Malleable Actions of Weakly Rigid Groups, part I}, preprint, math.OA/0305306.

\item{[P3]} S. Popa: {\it Some properties of the
symmetric enveloping algebras with applications to amenability and
property} T, Documenta Math. {\bf 4} (1999), 665-744.

\enddocument